\documentclass[12pt]{article}
\usepackage{amsmath}
\usepackage{amsbsy}
\usepackage{amsfonts}
\usepackage{mathtools}
\usepackage{bbm}
\usepackage{amsmath}
\usepackage{amstext}
\usepackage{amsthm}
\usepackage{bbold}
\usepackage{amssymb}
\usepackage[dvips]{graphicx}
\usepackage{amsfonts}
\usepackage{caption2}
\usepackage{amsfonts,mathrsfs}
\usepackage{makeidx}
\usepackage{fancyhdr}
\usepackage{epic,eepic,epsfig}
\usepackage{mdwlist}
\usepackage{amscd}
\usepackage{amsbsy}
\usepackage{epsfig,graphicx}
\usepackage{epsfig,mathrsfs,wasysym,stmaryrd}
\usepackage[all]{xypic}
\usepackage[knot,poly]{xy}
\usepackage{latexsym}
\usepackage{comment}

\theoremstyle{plain}

\newtheorem{theorem}{Theorem}

\newtheorem{lemma}{Lemma}

\newenvironment{pf}{\medskip\noindent{Proof:}
  \hspace{-.5cm}      \enspace}{\hfill \qed \newline \smallskip}

\date{}

\begin{document}

\begin{center}
\textbf{\LARGE{The number of rooted forests in circulant graphs}}
\vspace{12pt}

{\large\textbf{L.~A.~Grunwald,}}\footnote{{\small\em Sobolev Institute of Mathematics,
Novosibirsk State University, lfb\_o@yahoo.co.uk}}
{\large\textbf{I.~A.~Mednykh,}}\footnote{{\small\em Sobolev Institute of Mathematics,
Novosibirsk State University, ilyamednykh@mail.ru}}
\end{center}

\section*{Abstract}

In this paper, we develop a new method to produce explicit formulas for the number $f_{G}(n)$ of rooted spanning forests in the circulant graphs $ G=C_{n}(s_1,s_2,\ldots,s_k)$ and $ G=C_{2n}(s_1,s_2,\ldots,s_k,n).$ These formulas are expressed through Chebyshev polynomials. We prove that in both cases the number of rooted spanning forests can be represented in the form $f_{G}(n)=p\,a(n)^2,$ where $a(n)$ is an integer sequence and $p$ is a prescribed natural number depending on the parity of $n$.

Finally, we find an asymptotic formula for $f_{G}(n)$ through the Mahler measure of the associated Laurent polynomial $P(z)=2k+1-\sum\limits_{i=1}^k(z^{s_i}+z^{-s_i}).$
\bigskip

\noindent
\textbf{Key Words:} rooted tree, spanning forest, circulant graph, Laplacian matrix, Chebyshev polynomial, Mahler measure\\

\textbf{AMS classification:} 05C30, 39A12\\

\section{Introduction}

The famous Kirchhoff's Matrix Tree Theorem~\cite{Kir47} states that the number of spanning trees in a graph can be expressed as the product of its non-zero Laplacian eigenvalues divided by the number of vertices. Since then, a lot of papers devoted to the complexity of various classes of graphs were published. In particular, explicit formulae were derived for complete multipartite graphs~\cite{Cay89, Austin60}, almost complete graphs~\cite{Wein58}, wheels~\cite{BoePro}, fans~\cite{Hil74}, prisms~\cite{BB87}, ladders~\cite{Sed69}, M\"obius ladders~\cite{Sed70}, lattices \cite{SW00}, anti-prisms \cite{SWZ16}, complete prisms~\cite{Sch74} and for many other families. For the circulant graphs some explicit and recursive formulae are given in \cite{XiebinLinZhang,   ZhangYongGol, ZhangYongGolin}.

Along with the number of spanning trees in a given graph one can be interested in the number of rooted spanning forests in the graph. According to the classical result \cite{KelChel} (see also more recent papers \cite{ChebotSham}, \cite{Knill}) this value can be found with the use of determinant of the matrix $\det(I+L).$ Here $L$ is the Laplacian matrix of the graph. This invariant was calculated for various families of graphs. At the same time it is known very little about analytic formulas for the number of spanning forests. One of the first results was obtained by O. Knill \cite{Knill} who proved that the number of rooted spanning forests in the complete graph $K_n$ on $n$ vertices is equal to $(n+1)^{n-1}.$  The rooted spanning forests in bipatite graphs were enumerated in \cite{YinJ_ChuL}. Explicit formulas for the number of rooted spanning forests for cyclic, star, line and some others graphs were given by \cite{Knill}. As for the number of unrooted forests,  it has much more complicated structure \cite{Calan, LiuChow, LajTak}.

Starting with Boesch and Prodinger \cite{BoePro} the idea to apply Chebyshev polynomials for counting various invariants of graphs arose. This idea provided a way to find complexity of circulant graphs and their natural generalisations in~\cite{KwonMedMed, Med1, MedMed2018,  ZhangYongGolin}.

Recently, asymptotical behavior of complexity for some families of graphs was investigated from the point of view of so called Malher measure \cite{GutRog, SilWil1}. Mahler measure of a polynomial $P(z),$ with complex coefficients, is the absolute value of the product of all roots of $P(z)$ whose modulus is greater than $1$ multiplied by the leading coefficient. For general properties of the Mahler measure see the paper \cite{Smyth08}.

The purpose of this paper is to present new formulas for the number of rooted spanning forests in circulant graphs and investigate their arithmetical properties and asymptotic.

We arrange the paper in the following way. First, in the sections~\ref{count} and \ref{oddcomplexity} we present new explicit formulas for the number of spanning forests in the undirected circulant graphs $C_{n}(s_1,s_2,\ldots,s_k)$ and $C_{2n}(s_1,s_2,\ldots,s_k,n)$ of even and odd valency respectively. They will be given in terms of Chebyshev polynomials. Next, in the section~\ref{circarithm} some arithmetic properties of the number of spanning forests are investigated. More precisely, it is shown that the number of spanning forests of the circulant graph $G$ can be represented in the form $f_{G}(n)=p\,a(n)^{2},$ where $a(n)$ is an integer sequence and $p$ is a prescribed natural number depending only of the parity of $n.$ At last, in the section~\ref{assection}, we use explicit formulas for $f_{G}(n)$ in order to produce its asymptotic in terms of Mahler measure of the associated polynomials. For circulant graphs of even valency the associated polynomial is $P(z)=2k+1-\sum\limits_{j=1}^{k}(z^{s_j}+z^{-s_j}).$ In this case  (Theorem~\ref{asymptotic1}) we have $f_{G}(n)\sim A^n,\,n\to\infty,$ where $A=M(P)$ is the Mahler measure of $P(z).$ For circulant graphs of odd valency we use the polynomial $R(z)=P(z)(P(z)+2).$ Then the respective asymptotics (Theorem~\ref{asymptotic3}) is $f_{G}(n)\sim K^{n},\,n\to\infty,$ where $K=M(R).$ In the last section~\ref{tables}, we illustrate the obtained results by a series of examples.

\section{Basic definitions and preliminary facts}

Consider a finite graph $G$ without loops. We denote the vertex and edge sets of $G$ by $V(G)$ and $E(G),$ respectively. Given $u,\,v\in V(G),$ we set $a_{uv}$ to be equal to the number of edges between vertices $u$ and $v.$ The matrix $A=A(G)=\{a_{uv}\}_{u,\,v\in V(G)}$ is called \textit{the adjacency matrix} of the graph $G.$ The degree $d(v)$ of a vertex $v\in V(G)$ is defined by $d(v)=\sum_{u\in V(G)}a_{uv}.$ Let $D=D(G)$ be the diagonal matrix with the elements $d_{vv}=d(v).$ The matrix $L=L(G)=D(G)-A(G)$ is called \textit{the Laplacian matrix}, or simply \textit{Laplacian}, of the graph $G.$

By $I_{n}$ we denote the identity matrix of order $n.$

Denote by $\chi_{G}(\lambda)=\det(\lambda I_{n}-L(G))$ the characteristic polynomial of the Laplacian matrix of the graph $G.$ Its extended form is $$\chi_{G}(\lambda)=\lambda^{n}+c_{n-1}\lambda^{n-1}+\ldots+c_{1}\lambda.$$
The theorem by Kelmans and Chelnokov \cite{KelChel} states that the absolute value of coefficient $c_{k}$ of $\chi_{G}(\lambda)$ coincides with the number of rooted spanning $k-$forests in the graph $G.$ Since all the Laplacian eigenvales of $G$ are non-negative, one can note that the sequence $c_{k}$ is alternating. So, the number of rooted spanning forests of the graph $G$ can be found by the formula
\begin{eqnarray}\label{KelmCheln}
f_{G}(n)&=&f_{1}+f_{2}+\ldots+f_{n}=|c_{1}-c_{2}+c_{3}-\ldots+(-1)^{n-1}|\\
\nonumber&=&(-1)^{n}\chi_{G}(-1)=\det(I_{n}+L(G)).
\end{eqnarray}
This result was independently obtained by many authors (P. Chebatorev and E. Shamis \cite{ChebotSham}, O. Knill \cite{Knill} and others).

Let $s_1,s_2,\ldots,s_k$ be integers such that $1\leq s_1<s_2<\ldots<s_k\leq\frac{n}{2}.$ The graph $C_{n}(s_1,s_2,\ldots,s_k)$ with $n$ vertices $0,1,2,\ldots,~{n-1}$ is called \textit{circulant graph} if the vertex $i,\, 0\leq i\leq n-1$ is adjacent to the vertices $i\pm s_1,i\pm s_2,\ldots,i\pm s_k\,(\textrm{mod}\ n).$ When $s_k<\frac{n}{2}$ all vertices of a graph have even degree $2k.$ If $n$ is even and $s_k=\frac{n}{2},$ then all vertices have odd degree $2k-1.$

We call an $n\times n$ matrix \textit{circulant,} and denote it by $circ(a_0, a_1,\ldots,a_{n-1})$ if it is of the form
$$circ(a_0, a_1,\ldots, a_{n-1})=
\left(\begin{array}{ccccc}
a_0 & a_1 & a_2 & \ldots & a_{n-1} \\
a_{n-1} & a_0 & a_1 & \ldots & a_{n-2} \\
  & \vdots &   & \ddots & \vdots \\
a_1 & a_2 & a_3 & \ldots & a_0\\
\end{array}\right).$$

It easy to see that adjacency and Laplacian matrices of the circulant graph are circulant matrices. The converse is also true. If the Laplacian matrix of a graph is circulant then the graph is also circulant.

Recall \cite{PJDav} that the eigenvalues of matrix $C=circ(a_0,a_1,\ldots,a_{n-1})$ are given by the following simple formulas $\lambda_j=P(\varepsilon^j_n),\,j=0,1,\ldots,n-1,$ where $P(x)=a_0+a_1 x+\ldots+a_{n-1}x^{n-1}$ and $\varepsilon_n$ is an order $n$ primitive
root of the unity. Moreover, the circulant matrix $T=circ(0,1,0,\ldots,0)$ is the matrix representation of the shift operator
$T:(x_0,x_1,\ldots,x_{n-2},x_{n-1})\rightarrow(x_1, x_2,\ldots,x_{n-1},x_0).$

Let $P(z)=a_{0}z^{d}+\ldots+a_{d}=a_{0}\prod\limits_{i=1}^d(z-\alpha_i)$ be a nonconstant polynomial with complex coefficients. Then, following Mahler \cite{Mahl62} its {\it Mahler measure} is defined to be
\begin{equation} M(P):=\exp(\int_0^1\log|P(e^{2\pi i t})|dt).\end{equation}
The value $M (P)$ had appeared earlier in a paper by Lehmer \cite{Lehm33}, in an alternative form
\begin{equation} M(P)=|a_0|\prod\limits_{|\alpha_i|>1}|\alpha_i|.\end{equation}

The concept of Mahler measure can be naturally extended to the class of Laurent polynomials $P(z)=a_{0}z^{p}+a_{1}z^{p+1}+\ldots+a_{s-1}z^{p+s-1}+a_{s}z^{p+s}= a_{s}z^p\prod\limits_{i=1}^{s}(z-\alpha_i),$ where $a_s\neq0$ and $p$ is an arbitrary, but not necessarily positive integer.

Let $T_{n}(z)=\cos(n\arccos z)$ be the Chebyshev polynomial of the first kind. We will use the following property of the Chebyshev polynomials $T_{n}(\frac{1}{2}(z+z^{-1}))=\frac{1}{2}(z^{n}+z^{-n}).$ See \cite{MasHand} for more general properties.

\section{The number of rooted spanning forests in circulant graphs of even valency}\label{count}

The aim of this section is to find new formulas for the numbers of rooted spanning forests of circulant graph $C_{n}(s_1,s_2,\ldots,s_k)$ in terms of Chebyshev polynomials. Here and below, we will use $G$ to denote the circulant graph under consideration.

\begin{theorem}\label{theorem1}
The number of rooted spanning forests $f_{G}(n)$ in the circulant graph $G=C_{n}(s_1,s_2,\ldots,s_k),$ $1\le s_1< s_2<\ldots< s_k<\frac{n}{2},$ is given by the formula
$$f_{G}(n)=\prod_{p=1}^{s_k}|2T_{n}(w_{p})-2|,$$ thereby $w_{p},\,p=1,2,\ldots,s_{k}$ are all the roots of the algebraic equation $\sum\limits_{j=1}^{k}(2T_{s_{j}}(w)-2)=1,$ where and $T_{s}(w)$ is the Chebyshev polynomial of the first kind.
\end{theorem}

\medskip\noindent\textbf{Proof:} The number of rooted spanning forests of the graph $G$ can be found by the formula $f_{G}=\det(I_{n}+L(G)).$ The latter value is equal to the product of all eigenvalues of the matrix $I_{n}+L(G).$ We denote by $T=circ(0,1,\ldots,0)$ the $n \times n$ cyclic shift operator. Consider the Laurent polynomial $P(z)=2k+1-\sum\limits_{i=1}^k(z^{s_i}+z^{-s_i}).$ Then the matrix $I_{n}+L(G)$ has the following form
$$I_{n}+L(T)=P(T)=(2k+1)I_n-\sum\limits_{i=1}^k(T^{s_i}+T^{-s_i}).$$
The eigenvalues of circulant matrix $T$ are $\varepsilon_n^j,\,j=0,1,\ldots,n-1,$ where $\varepsilon_n=e^\frac{2\pi i}{n}.$ Since all of them are distinct, the matrix $T$ is conjugate to the diagonal matrix $\mathbb{T}=diag(1,\varepsilon_n,\ldots,\varepsilon_n^{n-1})$ with diagonal entries $1,\varepsilon_n,\ldots,\varepsilon_n^{n-1}.$ So the matrix $I_{n}+L(G)$ is conjugate to the diagonal matrix $P(\mathbb{T}).$ This essentially simplifies the problem of finding eigenvalues of $I_{n}+L(G).$ Indeed, let $\lambda$ be an eigenvalue of $I_{n}+L(G)$ and $x$ be the respective eigenvector. Then we have the following system of linear equations
$$((2k+1-\lambda)I_n-\sum\limits_{i=1}^k(\mathbb{T}^{s_i}+\mathbb{T}^{-s_i}))x=0.$$ Recall the matrices under consideration are diagonal and the $(j+1,j+1)$-th entry of $\mathbb{T}$ is equal to $\varepsilon_{n}^{j},$ where $\varepsilon_{n}=e^{\frac{2\pi i}{n}}.$ Then, for any $j=0,\ldots, n-1,$ matrix $P(\mathbb{T})$ has an eigenvalue
$\lambda_{j}=P(\varepsilon_{n}^{j})=2k+1-\sum\limits_{i=1}^k(\varepsilon_{n}^{j s_i}+\varepsilon_{n}^{-j s_{i}}).$ Hence we have
\begin{equation}\label{rospforfor}
f_{G}(n)=\prod\limits_{j=0}^{n-1}P(\varepsilon_{n}^{j}).
\end{equation}

To continue the proof of the theorem we need the following lemma.
\begin{lemma}\label{result} We have
$$\prod\limits_{j=0}^{n-1}P(\varepsilon_{n}^{j}) =\prod\limits_{p=1}^{s_{k}}|2T_{n}(w_{p})-2|,$$
where $w_{p},\,j=1,\ldots,s_{k}$ are all the roots of the algebraic equation $\sum\limits_{j=1}^{k}(2T_{s_{j}}(w)-2)=1.$
\end{lemma}

To prove the above formula we use some statements from theory of resultants. We introduce integer polynomial $\widetilde{P}(z)=-z^{s_{k}}P(z).$ We note that $\widetilde{P}(z)$ is a monic polynomial with the same roots as $P(z)$ and its degree is $2s_{k}.$ As $P(z)=P(\frac{1}{z}),$ the roots look like $z_{1},\,\frac{1}{z_{1}},\ldots,z_{s_{k}},\,\frac{1}{z_{s_{k}}}.$

We have $\prod\limits_{j=0}^{n-1}P(\varepsilon_{n}^{j})= \prod\limits_{j=0}^{n-1}(-\varepsilon_{n}^{-s_{k}j}\widetilde{P}(\varepsilon_{n}^{j})) =(-1)^{(s_k+1)(n+1)-1}\prod\limits_{j=0}^{n-1}\widetilde{P}(\varepsilon_{n}^{j}).$  By the basic properties of resultants

\begin{eqnarray*}\prod\limits_{j=0}^{n-1}\widetilde{P}(\varepsilon_{n}^{j})&=&\textrm{Res}\,(\widetilde{P}(z),\,z^{n}-1)=\textrm{Res}\,(z^{n}-1,\,\widetilde{P}(z))\\
&=& \prod\limits_{z:\widetilde{P}(z)=0}(z^{n}-1)=\prod\limits_{z:P(z)=0}(z^{n}-1)\\
&=&\prod\limits_{p=1}^{s_{k}}(z_{p}^{n}-1)(z_{p}^{-n}-1) =(-1)^{s_{k}}\prod\limits_{p=1}^{s_{k}}(2T_{n}(w_{p})-2).
\end{eqnarray*}

We use the identity $T_{n}(\frac{1}{2}(z+z^{-1}))=\frac{1}{2}(z^{n}+z^{-n}).$ Here $w_{p}=\frac{1}{2}(z_{p}+\frac{1}{z_{p}}),\,p=1,\ldots,s_{k}.$ These numbers are the roots of the algebraic equation $\sum\limits_{j=1}^{k}(2T_{s_j}(w)-2)=1.$  Since the righthand side of equation~(\ref{rospforfor}) is a positive integer, the lemma is proved.

By making use of Lemma~\ref{result}, we finish the proof of the theorem.
$\hfill \qed$





\section{The number of rooted spanning forests in circulant graphs of odd valency}\label{oddcomplexity}

This section is devoted to investigation of the numbers of rooted spanning forests
in circulant graph $C_{2n}(s_{1},s_{2},\ldots,s_{k},n)$ in terms of Chebyshev polynomials.


\begin{theorem}\label{odddegree}
Let $G=C_{2n}(s_{1},s_{2},\ldots,s_{k},n),\,1\leq s_{1}<s_{2}<\ldots<s_{k}<n,$ be a circulant graph of odd degree. Then the number $f_{G}(n)$ of rooted spanning forests in the graph $G$ is given by the formula
$$f_{G}(n)=\prod_{p=1}^{s_k}(2T_n(u_p)-2)(2T_n(v_p)+2),$$
where the numbers $u_p\text{ and } v_p,\,p=1,2,\ldots,s_k$ are respectively the roots of the algebraic equations $Q(u)-1=0$ and $Q(v)+1=0,$ where $Q(w)=2k+2-2\sum\limits_{i=1}^{k}T_{s_i}(w)$ and $T_k(w)$ is the Chebyshev polynomial of the first kind.
\end{theorem}

\begin{pf} In order to find the number of rooted spanning forests $f_{G}(n)$ in the graph $C_{2n}(s_{1},s_{2},\ldots,s_{k},n)$ we need to evaluate the determinant $\det(I_{2n}+L(G)).$ One can be represented the matrix $I_{2n}+L(G)$ in the form $$I_{2n}+L(G)=(2k+2)I_{2n}-\sum\limits_{j=1}^{k}(T^{s_j}+T^{-s_j})-T^n,$$ where $T$ is $2n\times 2n$ circulant matrix $circ(0,1,0,\ldots,0).$ The eigenvalues of circulant matrix $T$ are $\varepsilon_{2n}^j,\,j=0,1,\ldots,2n-1,$ where $\varepsilon_{2n}=e^\frac{2\pi i}{2n}.$ Since all of them are distinct, the matrix $T$ is conjugate to the diagonal matrix $\mathbb{T}=diag(1,\varepsilon_{2n},\ldots,\varepsilon_{2n}^{2n-1})$ with diagonal entries $1,\varepsilon_{2n},\ldots,\varepsilon_{2n}^{2n-1}$. To find the determinant $\det(I_{2n}+L(G))$ we use the product of all eigenvalues of matrix $I_{2n}+L(G).$ The matrix $I_{2n}+L(G)$ is conjugate to the diagonal matrix with eigenvalues
$$\lambda_j=2k+2-\sum\limits_{l=1}^{k}(\varepsilon_{2n}^{j\,s_l}+
\varepsilon_{2n}^{-j \,s_l})-\varepsilon_{2n}^{j n},\,j=0,1,\ldots,2n-1.$$
All of them are non-zero.

Consider the following Laurent polynomial $P(z)=2k+2-\sum\limits_{i=1}^{k}(z^{s_i}+z^{-s_i}).$ Since $\varepsilon_{2n}^n=-1,$  we can write $\lambda_j=P(\varepsilon_{2n}^{j})-1$ if $j$ is even and $\lambda_j=P(\varepsilon_{2n}^{j})+1$ if $j$ is odd. By the formula~\ref{KelmCheln} we have
$$f_{G}(n)=\prod\limits_{j=0}^{2n-1}\lambda_{j}=\prod\limits_{s=0}^{n-1}
(P(\varepsilon_{2n}^{2s})-1)\prod\limits_{s=0}^{n-1}(P(\varepsilon_{2n}^{2s+1})+1)$$
$$=\prod\limits_{s=0}^{n-1}(P(\varepsilon_{2n}^{2s})-1)\frac{\prod\limits_{p=0}^{2n-1}
(P(\varepsilon_{2n}^{p})+1)}{\prod\limits_{s=0}^{n-1}(P(\varepsilon_{2n}^{2s})+1)}=
\prod\limits_{s=0}^{n-1}(P(\varepsilon_{n}^{s})-1)\frac{\prod\limits_{p=0}^{2n-1}
(P(\varepsilon_{2n}^{p})+1)}{\prod\limits_{s=0}^{n-1}(P(\varepsilon_{n}^{s})+1)}.$$
\end{pf}

By making use of Lemma~\ref{result} and arguments from the proof of Theorem \ref{theorem1} we obtain
\begin{enumerate}
\item[(i)]
$\prod\limits_{s=0}^{n-1}(P(\varepsilon_{n}^{s})-1)=(-1)^{n(s_{k}+1)}
\prod_{p=1}^{s_k}(2T_n(u_p)-2),$

\item[(ii)]
$\prod\limits_{s=0}^{n-1}(P(\varepsilon_{n}^{s})+1)=(-1)^{n(s_{k}+1)}
\prod_{p=1}^{s_k}(2T_n(v_p)-2),$  and

\item[(iii)]
$\prod\limits_{p=0}^{2n-1}(P(\varepsilon_{2n}^{p})+1)=\prod_{p=1}^{s_k}(2T_{2n}(v_p)-2),$
\end{enumerate}
where $u_p$ and $v_p$ are the same as in the statement of the theorem. Hence,
$$f_{G}(n)=\prod_{p=1}^{s_k}(2T_n(u_p)-2)
\prod_{p=1}^{s_k}\frac{T_{2n}(v_p)-1}{T_n(v_p)-1}.$$
Finally, taking into account the identity
$T_{2n}(w)-1=2(T_n(w)-1)(T_n(w)+1)$ we get
$$f_{G}(n)=\prod_{p=1}^{s_{k}}(2T_{n}(u_{p})-2)(2T_{n}(v_{p})+2).$$
\hfill\qed

\section{Arithmetic properties of the number of rooted spanning forests for circulant graphs}\label{circarithm}

It has been proved in the paper \cite{MedMed2018} that the number of spanning trees $\tau(n)$ in circulant graph $C_{n}(s_{1}, s_{2},\dots, s_{k})$ is given by the formula $\tau(n)=p\,n\,a(n)^2,$ where $a(n)$ is an integer sequence and $p$ is a prescribed natural number depending only of pairity of $n$. The aim of the next theorem is to find a similar phenomenon for the number of rooted spanning forests.

Recall that any positive integer $p$ can be uniquely represented in the form $p=q \,r^2,$ where $p$ and $q$ are positive integers and $q$ is square-free. We will call $q$ the \textit{square-free part} of $p.$

\begin{theorem}\label{lorenzini}
Let $f_{G}(n)$ be the number of spanning forests in the circulant graph $$C_{n}(s_1,s_2,\ldots,s_k),\,1\le s_1<s_2<\ldots<s_k<\frac{n}{2}.$$ Denote by $p$ the number of odd elements in the sequence $s_1,s_2,\ldots,s_k$ and let $q$ be the square-free part of $4p+1.$ Then there exists an integer sequence $a(n)$ such that
\begin{enumerate}
\item[$1^0$] $f_{G}(n)=a(n)^2,$ if $n$ is odd;
\item[$2^0$] $f_{G}(n)=q\,a(n)^2,$ if $n$ is even.
\end{enumerate}
\end{theorem}
\begin{pf} The number of odd elements in the sequence $s_1,s_2,s_3,\ldots,s_k$ is counted by the formula $p=\sum\limits_{i=1}^{k}\frac{1-(-1)^{s_i}}{2}.$

We already know that all eigenvalues of the $I_{n}+L(G)$ are given by the formulas $\lambda_j=P(\varepsilon_n^{j}),\,j=0,\ldots,n-1,$ where $P(z)=2k+1-\sum\limits_{i=1}^{k}(z^{s_i}+z^{-s_i})$ and $\varepsilon_n=e^{\frac{2\pi \texttt{i}}{n}}.$ We note that $\lambda_{n-j}=P(\varepsilon_n^{n-j})=P(\varepsilon_n^{j})=\lambda_j.$

Since $\lambda_{0}=P(\varepsilon_{n}^{0})=P(1)=1,$ by the formula~(\ref{KelmCheln}) we have $f_{G}=\prod\limits_{j=1}^{n-1}\lambda_j.$ Since $\lambda_{n-j}=\lambda_j,$  we obtain $f_{G}=(\prod\limits_{j=1}^{\frac{n-1}{2}}\lambda_j)^2$ if $n$ is odd and $f_{G}=\lambda_{\frac{n}{2}}(\prod\limits_{j=1}^{\frac{n}{2}-1}\lambda_j)^2$ if $n$ is even. We note that each algebraic number $\lambda_j$ comes with all its Galois conjugate \cite{Lor}. So, the numbers $c(n)=\prod\limits_{j=1}^{\frac{n-1}{2}}\lambda_j$ and $d(n)=\prod\limits_{j=1}^{\frac{n}{2}-1}\lambda_j$ are integers. Also, for even $n$ we have $\lambda_{\frac{n}{2}}=P(-1)=2k+1-\sum\limits_{i=1}^{k}((-1)^{s_i}+(-1)^{-s_i}) =1+2\sum\limits_{i=1}^{k}(1-(-1)^{s_i})=4p+1.$ Hence, $f_{G}= c(n)^2$ if $n$ is odd and $f_{G}=(4p+1)\,d(n)^2$ if $n$ is even. Let $q$ be the free square part of $4p+1$ and $4p+1=q\,r^2.$ Setting $a(n)=c(n)$ in the first case and $a(n)=r\,d(n)$ in the second, we conclude that number $a(n)$ is always integer and the statement of theorem follows.
\end{pf}
\bigskip

The following theorem clarifies some number-theoretical properties of the number of rooted spanning forest $f_{G}(n)$ for circulant graphs of odd valency.
\bigskip

\begin{theorem}\label{lorenzinin} Let $f_{G}(n)$ be the number of rooted spanning forests in the circulant graph
$$G=C_{2n}(s_1,s_2,s_3,\ldots,s_k,n),\,1\le s_1< s_2<\ldots< s_k<n.$$ Denote by $p$ the number of odd elements in the sequence $s_1,s_2,s_3,\ldots,s_k.$ Let $q$ and $r$ be the square-free part of $4p+1$ and $4p+3$ respectively. Then there exists an integer sequence $a(n)$ such that
\begin{enumerate}
\item[$1^{0}.$] $f_{G}(n)=q\,a(n)^2,$ if $n$ is even;
\item[$2^{0}.$] $f_{G}(n)=r\,a(n)^2,$ if $n$ is odd.
\end{enumerate}
\end{theorem}
\begin{pf} The number $p$ of odd elements in the sequence $s_1,s_2,\ldots,s_k$ is equal to $\sum\limits_{i=1}^{k}\frac{1-(-1)^{s_i}}{2}.$ The eigenvalues of the matrix $I_{2n}+L(G)$ are given by the formulas $$\lambda_{j}=1+P(\varepsilon_{2n}^{j})-(-1)^{j},\,0=1,2,\ldots,2n-1,$$
where $P(z)=2k+1-\sum\limits_{l=1}^{k}(z^{s_l}+z^{-s_l})$ and $\varepsilon_{2n}=e^{\frac{\pi i}{n}}.$

Since $\lambda_{0}=1+P(1)-1=1$ by the formula~\ref{KelmCheln} we have $f_{G}(n)=\prod\limits_{j=1}^{2n-1}\lambda_{j}.$ Since $\lambda_{2n-j}=\lambda_j,$ we obtain $f_{G}(n)=\lambda_{n}(\prod\limits_{j=1}^{n-1}\lambda_{j})^2,$ where $\lambda_{n}=1+P(-1)-(-1)^n.$ Now we have
$$\lambda_{n}=2k+2-(-1)^n-2\sum\limits_{l=1}^{k}(-1)^{s_l} =2-(-1)^{n}+4\sum\limits_{l=1}^{k}\frac{1-(-1)^{s_{l}}}{2}=4\,p+2-(-1)^{n}.$$
So, $\lambda_{n}=4\,p+1,$ if $n$ is even and $\lambda_{n}=4\,p+3,$ if $n$ is odd. We note that each algebraic number $\lambda_{j}$ comes into the product  $\prod\limits_{j=1}^{n-1}\lambda_{j}$ together with all its Galois conjugate, so the number $c(n)=\prod\limits_{j=1}^{n-1}\lambda_{j}$ is an integer \cite{Lor}.

Hence, $f_{G}(n)=(4\,p+1)c(n)^2,$ if $n$ is even and $f_{G}(n)=(4\,p+3)\,c(n)^2,$ if $n$ is odd. Let $q$ and $r$ be the free square parts of $4\,p+1$ and of $4p+3$ respectively. Then for some integers $x$ and $y$ we have $4\,p+1=q\,x^2$ and $4\,p+3=r\,y^2.$

Now, the number $f_{G}(n)$ can be represented in the form
\begin{enumerate}
\item $\displaystyle{f_{G}(n)=q\,(x\,c(n))^{2}}$ if $n$ is even and
\item $\displaystyle{f_{G}(n)=r\,(y\,c(n))^{2}}$ if $n$ is odd.
\end{enumerate}
Setting $a(n)=x\,c(n)$ in the first case and $a(n)=y\,c(n)$ in the second, we conclude that number $a(n)$ is always integer. The theorem is proved.
\end{pf}

\section{Asymptotics for the number of spanning forests}\label{assection}

In this section we give asymptotic formulas for the number of rooted spanning forests in circulant graphs.


\begin{theorem}\label{asymptotic1} The number of rooted spanning forests in the circulant graph $G=C_{n}(s_1,s_2,\ldots,s_k),$ $1\le{s_1}<s_2<\ldots<s_k<\frac{n}{2}$
has the following asymptotics
$$f_{G}(n)\sim A^n,\text{ as }n\to\infty,$$
where $A=\exp(\int_0^1\log(P(e^{2\pi i t}))dt)$
is the Mahler measure of Laurent polynomial
$P(z)=2k+1-\sum\limits_{i=1}^k(z^{s_i}+z^{-{s_i}}).$

\end{theorem}
\begin{pf} By Theorem~\ref{theorem1} the number of rooted spanning forests $f_{G}(n)$ is given by
$$f_{G}(n)=\prod_{p=1}^{s_k}|2\,T_n(w_p)-2|.$$

We have $T_{n}(w_s)=\frac{1}{2}(z_{s}^{n}+z_{s}^{-n}),$ where the $z_{s}$ and $1/z_{s}$ are all the roots of the polynomial $P(z).$ If $\varphi\in\mathbb{R}$ then $P(e^{i\,\varphi})=2k+1-\sum\limits_{j=1}^{k} (e^{i s_{j}\,\varphi}+e^{-i s_{j}\,\varphi}) =2k+1-2\sum\limits_{j=1}^{k}\cos(s_{j}\,\varphi)\geq1,$ so  $|z_{s}|\neq1,\,s=1,2,\ldots,s_{k}.$ Replacing $z_s$ by $1/z_s,$ if it is necessary, we can assume that $|z_s|>1$ for all $s=1,2,\ldots,s_k.$ Then $T_n(w_s)\sim\frac{1}{2}z_s^n,$ as $n$ tends to $\infty.$ So,
$|2T_n(w_s)-2|\sim|z_s|^n,\,\,n\to\infty.$
Hence
$$\prod_{s=1}^{s_k}|2\,T_n(w_s)-2|\sim \prod_{s=1}^{s_k}|z_s|^n=
\prod\limits_{P(z)=0,\,|z|>1}|z|^n=A^n,$$
where $A=\prod\limits_{P(z)=0,\,|z|>1}|z|$  is the Mahler measure of $P(z).$
By the results mentioned in the preliminary part, it can be found by formula
$A=\exp(\int_0^1\log(P(e^{2\pi i t}))dt).$

Finally, $$f_{G}(n)=\prod_{s=1}^{s_k}|2\,T_n(w_s)-2|\sim A^n,\,n\to\infty.$$
\end{pf}
\bigskip

The next theorem is a direct consequence of Theorem~\ref{odddegree} and can be proved by  the same arguments as Theorem~\ref{asymptotic1}.

\bigskip
\begin{theorem}\label{asymptotic3} The number of rooted spanning forests $f_{G}(n)$ in the circulant graph $G=C_{2n}(s_1,s_2,\ldots,s_k,n),\,1\le s_1< s_2<\ldots< s_k<n$ has the following asymptotic
$$f_{G}(n)\sim K^{n},\text{ as }n\to\infty.$$
Here $K=\exp(\int\limits_0^1\log(P^{2}(e^{2\pi i t})-1)dt)$ is the Mahler measure of the Laurent polynomial $(P(z)-1)(P(z)+1),$ where $P(z)=2k+2-\sum\limits_{i=1}^k(z^{s_i}+z^{-s_i}).$
\end{theorem}
\bigskip

\section{Examples}\label{tables}

\begin{enumerate}

\item[$1^{\circ}$] \textbf{Cycle graph $G=C_{n}(1)=C_{n}.$} We need to solve the equation $1+2-2T_{1}(w)=0.$ We have $w=3/2.$ So, $f_{G}(n)=2T_{n}(3/2)-2.$ Furthermore, $f_{G}(n)\underset{n\to\infty}{\sim}(\frac{3+\sqrt{5}}{2})^{n}.$ Also, we have $f_{G}(n)=5F_{n}^{2},$ if $n$ is even, and $f_{G}(n)=L_{n}^{2},$ if $n$ is odd, where $F_{n}$ and $L_{n}$ are the Fibonacci and Lucas numbers respectively.

\item[$2^{\circ}$] \textbf{Graph $G=C_{n}(1,2).$} We need to solve the equation  $1+4-2T_{1}(w)-2T_{2}(w)=0.$ Its roots are $w_{1}=\frac{1}{4}(-1+\sqrt{29})$ and $w_{2}=\frac{1}{4}(-1-\sqrt{29}).$

Hence, by Theorem~\ref{theorem1}, $f_{G}(n)=|2T_{n}(w_{1})-2|\cdot|2T_{n}(w_{2})-2|\underset{n\to\infty}{\sim}A^{n},$ where $A=\frac{1}{4}(7+\sqrt{5}+\sqrt{38+14\sqrt{5}})\simeq4.3902568\ldots.$ By Theorem~\ref{lorenzini}, there exists an integer sequence $a(n)$ such that $f_{G}(n)=5\,a(n)^{2},$ if $n$ is even, and $f_{G}(n)=a(n)^{2},$ if $n$ is odd.

\item[$3^{\circ}$]\textbf{Graph $G=C_{n}(1,3).$} Let $w_{1},\,w_{2}$ and $w_{3}$ be the roots of the equation $1+4-2T_{1}(w)-2T_{3}(w)=0.$  Then $f_{G}(n)=|2T_{n}(w_{1})-2||2T_{n}(w_{2})-2||2T_{n}(w_{3})-2|\underset{n\to\infty}{\sim}A^{n},$ where $A\simeq4.48461\ldots$ is a suitable root of the equation $z^{4}-4z^{3}-2z^{2}-z+1=0.$ By Theorem~\ref{lorenzini}, $f_{G}(n)=a(n)^{2}$ for some integer sequence $a(n).$

\item[$4^{\circ}$]{\textbf{Graph M\"obius ladder $G=C_{2n}(1,n).$}}
We have to solve the equations $3-2T_{1}(w)=0$ and $5-2T_{1}(w)=0.$ Their roots are $u_{1}=\frac{3}{2}$ and $v_{1}=\frac{5}{2}$ respectively. Then $f_{G}(n)=(2T_{n}(\frac{3}{2})-2)(2T_{n}(\frac{5}{2})+2)\underset{n\to\infty}{\sim}A^{n},$ where $A=\frac{1}{4}(3+\sqrt{5})(5+\sqrt{21})\simeq 12.5438\ldots.$ By Theorem~\ref{lorenzinin}, $f_{G}(n)=5a(n)^{2},$ if $n$ is even, and $f_{G}(n)=7a(n)^{2},$ if $n$ is odd for some integer sequence $a(n)$.

\end{enumerate}

\section*{ACKNOWLEDGMENTS}
This work was supported by by the Russian Foundation for Basic Research (projects 18-01-00420, and 18-501-51021).

\end{document}